\documentclass[11pt]{amsart}
\usepackage{amsmath, amsthm, latexsym, amssymb, amsfonts, epsfig, epsf}
\usepackage[dvips]{color}

\newenvironment{pf}{\proof[\proofname]}{\endproof}
\theoremstyle{plain}
\newtheorem{Th}{Theorem}[section]
\newtheorem{Cor}[Th]{Corollary}

\numberwithin{equation}{section} \theoremstyle{definition}
\newtheorem{Rem}[Th]{Remark}

\newtheorem{Ex}{Example}
\newtheorem{Prob}{Problem}
\newtheorem{Alg}{Algorithm}
\newtheorem{Def}[Th]{Definition}

\newcommand{\cal}[1]{\mathcal{#1}}
\newcommand{\C}{\mathbb C}
\newcommand{\Z}{\mathbb Z}
\newcommand{\R}{\mathbb R}
\newcommand{\Pp}{\mathbb P}
\newcommand{\G}{\Gamma}

\newcommand{\D}{\Delta}
\newcommand{\cA}{\cal A}

\newcommand{\cO}{\cal O}

\newcommand{\cR}{\cal R}

\newcommand{\p}{\partial}
\newcommand{\T}{{\mathbb T}}
\newcommand{\Jac}{{J_f^{\,\T}}}
\newcommand{\gres}{{\cR_f^{\T}}}

\newcommand{\Sig}{\Sigma}

\newcommand{\res}{\operatorname{res}}

\newcommand{\Vol}{\operatorname{Vol}}

\newcommand{\rs}[1]{Section~\ref{S:#1}}

\newcommand{\rpr}[1]{Problem~\ref{Pr:#1}}
\newcommand{\rr}[1]{Remark~\ref{R:#1}}

\newcommand{\re}[1]{(\ref{e:#1})}
\newcommand{\rc}[1]{Corollary~\ref{C:#1}}
\newcommand{\rt}[1] {Theorem~\ref{T:#1}}

\newcommand{\rf}[1]{Figure~\ref{F:#1}}

\begin{document}


\title{Global residues for sparse polynomial systems}
\author{Ivan Soprunov}

\address{Department of Mathematics and Statistics,
University of Massachusetts, Amherst, USA}
\email{isoprou@math.umass.edu}

\begin{abstract} We consider families of sparse Laurent polynomials
$f_1,\dots,f_n$ with a finite set of common zeroes $Z_f$ in the
torus $\T^n=(\C-\{0\})^n$. The global residue assigns to every
Laurent polynomial $g$ the sum of its Grothendieck residues over
$Z_f$. We present a new symbolic algorithm for computing the
global residue as a rational function of the coefficients of the
$f_i$ when the Newton polytopes of the $f_i$ are full-dimensional.
Our results have consequences in sparse polynomial interpolation
and lattice point enumeration in Minkowski sums of polytopes.
\end{abstract}

\maketitle

\section{Introduction}
Let $f_1=\cdots=f_n=0$ be a system of Laurent polynomial equations in $n$
variables whose Newton polytopes are $\D_1,\dots,\D_n$. Suppose the
solution set $Z_f$ in the algebraic torus $\T=(\C-\{0\})^n$ is finite.
This will be true if the coefficients of the $f_i$ are generic.

The {\it global residue} assigns to every polynomial $g$ the sum
over $Z_f$ of Grothendieck residues of the meromorphic form
$$\omega_g=\frac{g}{f_1\dots f_n}\,\frac{dt_1}{t_1}\wedge\dots\wedge\frac{dt_n}{t_n}.$$
It is a symmetric function of the solutions and hence depends
rationally on the coefficients of the system. Global residues were
studied by Tsikh, Gelfond and Khovanskii,  and later by Cattani,
Dickenstein and Sturmfels \cite{Ts, GKh2, CaD, CDS}. There are
numerous applications of global residues ranging from elimination
algorithms in commutative algebra \cite{CDS} to integral
representation formulae in complex analysis \cite{Ts}.

There have been several approaches to the problem of computing the
global residue explicitly. Gelfond and Khovanskii considered
systems with {\it generically positioned} Newton polytopes. The
latter means that for any linear function $\xi$, among the
corresponding extremal faces $\D_1^\xi,\dots,\D_n^\xi$ at least
one is a vertex. In this case Gelfond and Khovanskii found an
explicit formula for the global residue as a Laurent polynomial in
the coefficients of the system \cite{GKh2}. In general, the global
residue is not a Laurent polynomial and one should not expect a
closed formula for it.

Another approach was taken by Cattani and Dickenstein in the {\it
generalized unmixed} case when all $\D_i$ are dilates of a
single $n$-dimensional polytope (see \cite{CaD}). Their algorithm requires
computing a certain element in the homogeneous coordinate ring of
a toric variety whose {\it toric residue} equals 1. Then the {\it
Codimension Theorem} for toric varieties \cite{CD} allows global
residue computation using Gr\"obner bases techniques. This approach was
then extended by D'Andrea and Khetan to a slightly more general {\it ample}
case \cite{DK}. However, it
turned out to be a hard (and still open) problem to find an
element of toric residue~1 for an arbitrary collection of
polytopes.\footnote{For a combinatorial construction of such an
element that provides a solution to the problem for a wide class
of polytopes see~\cite{KS}.}

In this paper we present a new algorithm for computing the global
residue when $\D_1,\dots,\D_n$ are arbitrary $n$-dimensional
polytopes (see \rs{algorithm}). The proof of its correctness is in
the spirit of Cattani and Dickenstein's arguments, but avoids the
toric residue problem. It relies substantially on the Toric
Euler--Jacobi theorem due to Khovanskii~\cite{uspehi}. Also in our
algorithm we replace Gr\"obner bases computations with solving
a linear system. This gives an expression of the global residue as
a quotient of two determinants. The same idea was previously used
by D'Andrea and Khetan in \cite{DK}.

There is an intimate relation between the global residue and
polynomial interpolation. In the classical case this was
understood already by Kronecker. In \rs{sparseinterpol} we show
how the Toric Euler--Jacobi theorem gives sparse polynomial
interpolation. Another consequence of this theorem is a lower
bound on the number of interior lattice points in the Minkowski
sum of $n$ full-dimensional lattice polytopes in $\R^n$
(\rc{bound}). Finally, our results give rise to interesting
optimization problems about Minkowski sums which we discuss in
\rs{geometry}.

\subsection*{Acknowledgments} This work was inspired by
discussions with Askold Khovanskii when I visited him in Toronto.
I am grateful to Askold for the valuable input and his
hospitality. I thank Amit Khetan for numerous discussions on this
subject and Eduardo Cattani and David Cox for their helpful
comments.

\section{Global Residues in $\C^n$ and Kronecker interpolation}

We begin with classical results on polynomial interpolation that
show how global residues come into play.

Let $f_1,\dots, f_n\in\C[t_1,\dots,t_n]$ be polynomials of degrees
$d_1,\dots,d_n$, respectively. Let $\rho=\sum d_i-n$ denote the
{\it critical degree} for the $f_i$. We assume that the solution
set $Z_f\subset\C^n$ of $f_1=\dots=f_n=0$ consists of a finite
number of simple roots. The following is a classical problem on
polynomial interpolation.

\begin{Prob} Given $\phi: Z_f\to\C$ find a polynomial $g\in\C[t_1,\dots,t_n]$ of
degree at most $\rho$ such that $g(a)=\phi(a)$ for all $a\in Z_f$.
\end{Prob}

Kronecker \cite{K} suggested the following solution to this problem. Since each $f_i$ vanishes
at every $a\in Z_f$ we can write (non-uniquely):
\begin{equation}\label{e:1}
f_i=\sum_{j=1}^n g_{ij}(t_j-a_j),\quad a=(a_1,\dots,a_n),\ \ g_{ij}\in\C[t_1,\dots,t_n].
\end{equation}
The determinant $g_a$ of the polynomial matrix $[g_{ij}]$ is a
polynomial of degree at most~$\rho$. Notice that the value of
$g_a$ at $t=a$ equals the value of the Jacobian
$J_f=\det\left(\frac{\p f_i}{\p t_j}\right)$ at $t=a$, since
$\frac{\p f_i}{\p t_j}(a)=g_{ij}(a)$ from \re{1}. Also $g_a(a')=0$
for any $a'\in Z_f$, $a'\neq a$. Indeed, substituting $t=a'$ into
\re{1} we get
\begin{equation}
0=\sum_{j=1}^n g_{ij}(a')(a'_j-a_j),\nonumber
\end{equation}
which means that the non-zero vector $a'-a$ is in the kernel of
the matrix $[g_{ij}(a')]$, hence $\det[g_{ij}(a')]=0$. It remains
to put
\begin{equation}\label{e:2}
g=\sum_{a\in Z_f}\frac{\phi(a)}{J_f(a)}\,g_a.
\end{equation}

\begin{Rem} When choosing $g_{ij}$ in \re{1} one can assume that
$$g_{ij}=\frac{\p F_i}{\p t_j}+\text{lower degree terms},$$ where $F_i$
is the main homogeneous part of $f_i$. Then \re{2} becomes
\begin{equation}\label{e:3}
g=\bigg(\sum_{a\in Z_f}\frac{\phi(a)}{J_f(a)}\bigg)J_F+\text{\rm lower degree terms},
\end{equation}
where $J_F$ is the Jacobian of the $F_i$.
\end{Rem}


\begin{Def} Given $g\in\C[t_1,\dots,t_n]$
the sum of the local Grothendieck residues
$$\cR_f(g)=\sum_{a\in Z_f}\res_a\left(\frac{g}{f_1\dots f_n}\,dt_1\wedge\dots\wedge dt_n\right)$$
is called the {\it global residue} of $g$
for the system $f_1=\dots=f_n=0$. In the case of simple roots of the system we get
$$\cR_f(g)=\sum_{a\in Z_f}\frac{g(a)}{J_f(a)}.$$
\end{Def}

\begin{Th} ({The Euler-Jacobi theorem}) Let $f_1=\dots=f_n=0$ be a generic polynomial system with
$\deg(f_i)=d_i$. Then for any $h$ of degree less than $\rho=\sum d_i-n$
the global residue $\cR_f(h)$ is zero.
\end{Th}

\begin{pf} Consider the function $h:Z_f\to\C$, $a\mapsto h(a)$. According to the previous
discussion there is a polynomial of the form \re{3}
which takes the same values on $Z_f$ as $h$. In other words,
$$h\equiv\bigg(\sum_{a\in Z_f}\frac{h(a)}{J_f(a)}\bigg)J_F+\text{\rm lower degree terms}\quad
(\text{mod }\ I),$$
where we used that the ideal $I=\langle f_1,\dots, f_n\rangle$ is radical and the roots are simple.
Comparing the homogeneous parts of degree $\rho$ in this equation we see that either
$J_F\in I$, which is equivalent to the $F_i$ having a non-trivial common zero (does not
happen generically), or the coefficient of $J_F$ (the global residue of $h$) is zero.
\end{pf}

\section{Global Residues in $\T^n$ and sparse polynomial
interpolation}\label{S:sparseinterpol}

Now we will consider {\it sparse} polynomial systems and define
the global residue in this situation. The word ``sparse''
indicates that instead of fixing the degrees of the polynomials we
fix their {\it Newton polytopes}. The Newton polytope $\D(f)$ of a
(Laurent) polynomial $f$ is defined as the convex hull in $\Z^n$
of the exponent vectors of all the monomials appearing in $f$.
Note that the Newton polytope of a generic polynomial of degree
$d$ is a $d$-dilate of the standard $n$-simplex. Such polynomials
are usually called {\it dense}. In what follows when we say
generic sparse polynomial we will mean that its Newton polytope is
fixed and the coefficients are generic.

Let $f_1,\dots, f_n\in\C[t_1^{\pm 1},\dots,t_n^{\pm 1}]$ be
Laurent polynomials whose Newton polytopes $\D_1,\dots, \D_n$ are
$n$-dimensional. We will assume that the solution set
$Z_f\subset\T^n$ of the system $f_1=\dots=f_n=0$ is finite. Here
$\T^n$ denotes the $n$-dimensional algebraic torus $(\C-\{0\})^n$.

Similar to the affine case define the {\it
global residue} of a Laurent polynomial $g$ as the sum
of the local Grothendieck residues
$$\gres(g)=\sum_{a\in Z_f}
   \res_a\left(\frac{g}{f_1\dots f_n}\,\frac{dt_1}{t_1}\wedge\dots\wedge\frac{dt_n}{t_n}\right).$$
When the roots of the system are simple we obtain
$$\gres(g)=\sum_{a\in Z_f}\frac{g(a)}{\Jac(a)},$$
where $\Jac=\det\left(t_j\frac{\p f_i}{\p t_j}\right)$ is
the {\it toric Jacobian} of the polynomials $f_1,\dots,f_n$.

Notice that $\cR^\T_f$ is a linear function on the space of
Laurent polynomials and depends rationally on the coefficients of
the $f_i$, since it is symmetric in the roots $Z_f$. In
\rs{algorithm} we will give an algorithm for computing $\gres(g)$
as a rational function of the coefficients of the system.

The following theorem due to A.~Khovanskii is a far reaching
generalization of the Euler--Jacobi theorem.

\begin{Th}\cite{uspehi}\label{T:Kh}
Let $f_1=\dots=f_n=0$ be a generic system of Laurent polynomials with $n$-dimensional
Newton polytopes $\D_1,\dots, \D_n$. Let $\D=\D_1+\dots+\D_n$ be
the Minkowski sum. Then
\begin{enumerate}
\item (Toric Euler-Jacobi) for any Laurent polynomial $h$ whose
Newton polytope lies in the interior of $\D$ the global residue
$\gres(h)$ is zero.
\item (Inversion of Toric Euler-Jacobi) for
any $\phi: Z_f\to \C$ with $\sum_{a\in Z_f}\phi(a)=0$ there exists
a polynomial $h$ whose Newton polytope lies in the interior of
$\D$ such that $\phi(a)=h(a)/\Jac(a)$.
\end{enumerate}
\end{Th}

Let us denote by $S_{\D^\circ}$ the vector space of all Laurent
polynomials whose Newton polytope lies in the interior $\D^\circ$
of $\D$. We have a linear map
\begin{equation}\label{e:linearmap}
\cA:S_{\D^\circ}\to\C^{|Z_f|},\quad
h\mapsto\bigg(\frac{h(a)}{\Jac(a)},\ a\in Z_f\bigg).
\end{equation}
Then the above theorem says that the image of $\cA$ is the
hyperplane $\{\sum x_i=0\}$ in~$\C^{|Z_f|}$. By Bernstein's
theorem the number of solutions $|Z_f|$ is equal to the normalized
mixed volume $n!\,V(\D_1,\dots,\D_n)$ of the polytopes
\cite{Bern}. We thus obtain a lower bound on the number of
interior lattice points of Minkowski sums.

\begin{Cor}\label{C:bound} Let $\D_1,\dots, \D_n$ be $n$-dimensional
lattice polytopes in $\R^n$ and $\D$ their Minkowski sum. Then the
number of lattice points in the interior of $\D$ is at least
$n!\,V(\D_1,\dots,\D_n)-1$.
\end{Cor}

It would be interesting to give a direct geometric proof of this
inequality and determine all collections $\D_1,\dots,\D_n$ for
which the bound is sharp. In the unmixed case $\D_1=\dots=\D_n=\D$
the inequality becomes
\begin{equation}\label{e:unmixedbound}
(n\D)^\circ\cap\Z^n\geq n!\Vol_n(\D)-1
\end{equation}
and can be deduced from the Stanley's Positivity theorem for the
Ehrhart polynomial \cite{Stan}. Recently Batyrev and Nill
described all possible $\D$ which give equality in \re{unmixedbound} (see \cite{BN}).
Here is a mixed case example which shows that the bound in \rc{bound} is sharp.

\begin{Ex} Let $\G(m)$ denote the simplex defined as the convex hull in $\R^n$ of
$n+1$ points $\{0,e_1,\dots,e_{n-1},me_n\}$, where $e_i$ is the $i$-th standard basis vector.

Consider a collection of $n$ such simplices $\G(m_1),\dots,\G(m_n)$
with $m_1\leq\dots\leq m_n$. It is not hard to see that their
mixed volume equals $m_n$. For example, one can consider a generic
system with these Newton polytopes and eliminate all but the last
variable to obtain a univariate polynomial of degree $m_n$. The
number of solutions of such a system is $m_n$, which is the mixed
volume by Bernstein's theorem.

Also one can see that the number of interior
lattice points in $\G(m_1)+\dots+\G(m_n)$ is exactly $m_n-1$. (In fact, these lattice
points are precisely the points $(1,\dots,1,k)$ for $1\leq k< m_n$.)
\end{Ex}

\begin{Cor} (Sparse Polynomial Interpolation)\label{C:interpol}
Let $f_1=\dots=f_n=0$ be a generic system with $n$-dimensional
Newton polytopes $\D_1,\dots,\D_n$ and let $\D$ be their Minkowski
sum. Let $Z_f\subset\T^n$ denote the solution set of the system.
Then for any function $\phi:Z_f\to\C$ there is a polynomial $g$
with $\D(g)\subseteq\D$ such that $g(a)=\phi(a)$. Moreover, $g$
can be chosen to be of the form $g=h+c\Jac$ for some $h$ with
$\D(h)\subset\D^\circ$ and a constant $c$.
\end{Cor}

\begin{pf} Consider a new function $\psi:Z_f\to\C$ given by
$$\psi(a)=\frac{\phi(a)}{\Jac(a)}-\frac{c}{MV},\quad \text{where }\
c=\sum_{a\in Z_f}\frac{\phi(a)}{\Jac(a)},\ MV=n!\, V(\D_1,\dots,\D_n).$$
Then the sum of the values of $\psi$ over the points of $Z_f$ equals zero. Therefore
there exists $h\in S_{\D^\circ}$ such that $h(a)=\Jac(a)\psi(a)=\phi(a)-\frac{c}{MV}\Jac(a)$
for all $a\in Z_f$. It remains to put $g=h+\frac{c}{MV}\Jac$.
\end{pf}

\begin{Rem} \rt{Kh} is an instance of a more general
result. Let $f_1,\dots,f_k$ be generic Laurent polynomials with
$n$-dimensional Newton polytopes $\D_1,\dots,\D_k$, for $k\leq n$.
The set of their common zeroes defines an algebraic variety $Z_k$
in $\T^n$. There is a way to embed $\T^n$ into a projective toric
variety $X$ so that the algebraic closure of $f_i=0$ defines a
Cartier divisor $D_i$ on $X$ and the closure $\overline Z_k$ is a
complete intersection in $X$. In \cite{Kh2} Khovanskii described
the space of top degree holomorphic forms on $\overline Z_k$. The
special case $k=n$ corresponds to the space of all functions on
the finite set $\overline Z_n=Z_f$ whose description we gave in
\rt{Kh}. It follows from cohomology computation on complete
intersections. In particular, for $k=n$, there is an exact
sequence of global sections
\begin{equation}\label{e:seq1}
\dots\to \bigoplus_{i=1}^n H^0(X,\cO(D-D_i+K))\to H^0(X,\cO(D+K))\to \C^{|Z_f|}\to \C \to 0.
\end{equation}
Here $D=D_1+\dots+D_n$, $K$ the canonical divisor, and $\cO(L)$
the invertible sheaf corresponding to a divisor $L$. The first non-zero map in \re{seq1} (from the right)
is the trace map, the second map is the residue
map we considered in \re{linearmap}, and the third one is given by $(f_1,\dots,f_n)$.
\end{Rem}

\section{Some commutative algebra}

As before consider a system of Laurent polynomial equations
$f_1=\dots=f_n=0$ whose Newton polytopes $\D_1,\dots,\D_n$ are
full-dimensional. For generic coefficients the system will have a
finite number of simple roots $Z_f$ in the torus $\T^n$. We concentrate
on the following problem.

\begin{Prob} Given a Laurent polynomial $g$ compute the global
residue $\gres(g)$ as a rational function of the coefficients of
the $f_i$.
\end{Prob}

We postpone the algorithm to the next section and now formulate
our main tool for solving the problem.

\begin{Th}\label{T:Cox} Let $\D_0,\dots, \D_n$ be $n+1$ full-dimensional
lattice polytopes in $\R^n$ and assume that $\D_0$ contains the origin in its interior.
Put
\begin{equation}\label{e:notation}
\tilde\D=\D_0+\dots+\D_n\quad\text{and}\quad
\tilde\D_{(i)}=\D_0+\dots+\D_{i-1}+\D_{i+1}+\dots+\D_n.
\end{equation}
Then for generic polynomials $f_i$ with Newton polytopes $\D_i$, for
$1\leq i\leq n$, the linear map
\begin{equation}\label{e:resmap}
\bigoplus_{i=0}^nS_{\tilde{\D}_{(i)}^\circ}\oplus\C\to
S_{\tilde{\D}^\circ}, \quad (h_0,\dots,h_n, c)\mapsto
h_0+\sum_{i=1}^n h_if_i+c\Jac
\end{equation}
is surjective.
\end{Th}

\begin{pf} Let $g$ be in $S_{\tilde{\D}^\circ}$. By \rc{interpol}
there exists a polynomial $h_0$ supported in
$\D^\circ=\tilde{\D}_{(0)}^\circ$ and a constant $c$ such that
$g(a)=h_0(a)+c\Jac(a)$ for all $a\in Z_f$, i.e. the polynomial
$g-h_0-c\Jac\in S_{\tilde{\D}^\circ}$ vanishes on $Z_f$. Now the
statement follows from \rt{Noether} below.
\end{pf}

The following statement can be considered as the toric version of
the classical Noether theorem in $\Pp^n$ (see \cite{Ts},
section 20.2).
\begin{Th}\label{T:Noether} Let $f_1,\dots,f_n$ be generic Laurent polynomials with
$n$-dimensional Newton polytopes $\D_1,\dots,\D_n$. Let $h$ be a
Laurent polynomial vanishing on $Z_f$. Assume $\D(h)$ lies in the
interior of $\tilde\D=\D_0+\D_1+\dots+\D_n$ for some $n$-dimensional
polytope~$\D_0$. Then $h$ can be written in the form
$$h=h_1f_1+\dots+h_nf_n,\quad\text{with }\ \D(h_i)\subset\tilde\D_{(i)}^\circ,$$
where $\tilde\D_{(i)}$ as in \re{notation}.
\end{Th}

\begin{pf} First we note that the statement remains true when
$\D_0=\{0\}$, i.e. $h$ is supported in the interior of
$\D=\D_1+\dots+\D_n$. This follows from the exact sequence
\re{seq1}. Indeed, if one considers the toric variety $X$
associated with $\D$ then each $f_i$ defines a semiample divisor
$D_i$ on $X$ with polytope $\D_i$. It is well-known that for any
semiample divisor $L$
\begin{equation}\label{e:fact}
H^{n-\dim\D_L}(X,\cO(L+K))\cong S_{\D_L^\circ},
\end{equation}
where $\D_L$ is the polytope of $L$ (see, for example \cite{Kh1,F}).
Thus the first term in \re{seq1} is isomorphic to
$S_{\D_{(1)}^\circ}\oplus\dots\oplus S_{\D_{(n)}^\circ}$, where
$\D_{(i)}=\sum_{j\neq i}\D_j$. Since the sequence is exact and $h$
lies in the kernel of the second map we get the required
representation.

Now assume $\D_0$ is $n$-dimensional. Let $X$ be the toric variety
associated with $\tilde\D$. Let $f_0$ be any monomial supported in
$\D_0$. Then $f_0,\dots,f_n$ define $n+1$ semiample divisors
$D_0,\dots, D_n$ on $X$ whose polytopes are $\D_0,\dots,\D_n$.
Since $f_1,\dots,f_n$ are generic (and so all their common zeroes
lie in $\T^n$) and $f_0$ is a monomial, the divisors $D_0,\dots,
D_n$ have empty intersection in $X$. Then the following twisted
Koszul complex of sheaves on $X$ is exact (see \cite{CD, DK}):
$$0\to\cO(K)\to\bigoplus_{i=0}^n\cO(D_i+K)\to\dots\to\bigoplus_{i=0}^n\cO(\tilde D-D_i+K)\to\cO(\tilde D+K)\to 0,$$
where $\tilde D=D_0+\dots+D_n$ and $K$ the canonical divisor on $X$.
The first few term of the cohomology sequence are
$$\dots\to\bigoplus_{i=0}^nH^0(X,\cO(\tilde D-D_i+K))\to H^0(X,\cO(\tilde D+K))\to 0,$$
where the middle map is given by $(f_0,\dots,f_n)$. This, with the help
of \re{fact}, provides
$$h=h_0f_0+h_1f_1+\dots+h_nf_n,\quad\text{where }\ \D(h_i)\subset\tilde\D_{(i)}^\circ.$$
Notice that $h_0$ vanishes on $Z_f$ and is supported in the
interior of $\tilde\D_{(0)}=\D$. Therefore, there exist $h_i'$ such
that
$$h_0=h_1'f_1+\dots+h_n'f_n,\quad\text{with }\ \D(h_i')\subset\D_{(i)}^\circ$$
by the previous case. It remains to note that
$\D(h_i'f_0)\subset\D_{(i)}^\circ+\D_0=\tilde\D_{(i)}^\circ$.
\end{pf}

\begin{Rem} \rt{Cox} has interpretation in terms of the homogeneous
coordinate ring $S_X$ of the toric variety $X$ associated with
$\tilde\D$ (see \cite{Coxhom}). One can homogenize $f_0,\dots, f_n$
to get elements $F_0,\dots,F_n\in S_X$ of {\it big} and {\it nef}
degrees. According to the Codimension Theorem of Cox and
Dickenstein (see \cite{CD}) the codimension of~$I=\langle
F_0,\dots F_n\rangle$ in {\it critical} degree (corresponding to
the interior of $\tilde\D$) equals 1. Then \rt{Cox} says that the
homogenization of the Jacobian $\Jac$ to the critical degree
generates the critical degree part of the quotient $S_X/I$.
\end{Rem}

\section{Algorithm for computing the global residue in
$\T^n$}\label{S:algorithm}

Now we will present an algorithm for computing the global residue
$\gres(g)$ for any Laurent polynomial $g$ assuming that the Newton
polytopes of the system are full-dimensional.

\begin{Alg} Let $f_1=\dots=f_n=0$ be a system of Laurent polynomial equations with
$n$-dimensional Newton polytopes $\D_1,\dots,\D_n$. As before we let $\D$ denote the Minkowski
sum of the polytopes.\newline
{\it Input:} A Laurent polynomial $g$ with Newton polytope $\D(g)$.
\begin{enumerate}
\item[{\it Step 1:}] Choose an $n$-dimensional polytope $\D_0$
with $0\in\D_0^\circ$ such that the Minkowski sum
$\tilde{\D}=\D_0+\D$ contains $\D(g)$ in its interior. \item[{\it
Step 2:}] Solve the system of linear equations
\begin{equation}\label{e:linsys}
g=h_0+\sum_{i=1}^n h_if_i+c\Jac
\end{equation}
for $c$, where $h_i$ are polynomials with unknown coefficients
supported in the interior of $\tilde\D_{(i)}$ (see \re{notation}).
\end{enumerate}
{\it Output:} The global residue $\gres(g)=c\,n!\,V(\D_1,\dots,\D_n)$.

\end{Alg}

\begin{pf} According to \rt{Cox}, given $g$ with $\D(g)\subset\tilde\D^\circ$
there exist $h_i$ supported in $\tilde\D_{(i)}^\circ$ and $c\in\C$ such that
$g=h_0+\sum_{i=1}^n h_if_i+c\Jac$. Taking the global residue we have
$$\gres(g)=\gres(h_0)+\gres\Big(\sum_{i=1}^n h_if_i\Big)+c\,\gres(\Jac)=c\,n!\,V(\D_1,\dots,\D_n),$$
where the first two terms vanish by \rt{Kh} (1) and the definition of the global residue.
\end{pf}

\begin{Rem}\label{R:jac} Notice that we can ignore those terms of $\Jac$ whose
exponents lie in the interior of $\D$ since their residue is zero by \rt{Kh} (1), and
work with the ``restriction'' of $\Jac$ to the boundary of $\D$.
\end{Rem}


We illustrate the algorithm with a small 2-dimensional example.

\begin{Ex} Consider a system of two equations in two unknowns.
$$\begin{cases}&\hspace{-10pt}f_1=a_1x+a_2y+a_3x^2y^2,\\ &\hspace{-10pt}f_2=b_1x+b_2xy^2+b_3x^2y^2.\end{cases}$$
The Newton polytopes $\D_1$, $\D_2$ and their Minkowski sum $\D$
are depicted in \rf{ex1}.

\begin{figure}[h]
\centerline{ \scalebox{0.80}{
\input{newex2.pstex_t}}}
\caption{}
\label{F:ex1}
\end{figure}

We compute the global residue of $g=x^5y^4$. Let $\D_0$ be the
triangle with vertices $(-1,0)$, $(0,-1)$ and $(2,1)$. Then the
Minkowski sum $\tilde\D=\D+\D_0$ contains $\D(g)=(5,4)$ in the
interior (see \rf{ex2}).
\begin{figure}[h]
\centerline{ \scalebox{0.75}{
\input{newex3.pstex_t}}}
\caption{}
\label{F:ex2}
\end{figure}
The vector space $S_{\tilde\D^\circ}$ has dimension 15 and a
monomial basis
\begin{equation}
S_{\tilde\D^\circ}=\langle xy, xy^2, xy^3, x^2, x^2y, x^2y^2,
x^2y^3, x^3y, x^3y^2, x^3y^3, x^3y^4, x^4y^2, x^4y^3, x^4y^4,
x^5y^4\rangle.\nonumber
\end{equation}
Now $\tilde\D_{(0)}=\D_1+\D_2=\D$, $\tilde\D_{(1)}=\D_0+\D_2$ and $\tilde\D_{(2)}=\D_0+\D_1$
and the corresponding vectors spaces are of dimension 4, 6 and 6, respectively. Here are their monomial bases:
$$ S_{\tilde{\D}_{(0)}^\circ}=\langle x^2y, x^2y^2, x^2y^3,
x^3y^3\rangle,\quad S_{\tilde{\D}_{(1)}^\circ}=\langle x, xy, xy^2,
x^2y, x^2y^2, x^3y^2\rangle, $$
$$S_{\tilde{\D}_{(2)}^\circ}=\langle
y, x, xy, x^2y, x^2y^2, x^3y^2\rangle,$$ as \rf{ex3} shows.
\begin{figure}[h]
\centerline{ \scalebox{0.75}{
\input{newex1.pstex_t}}}
\caption{}
\label{F:ex3}
\end{figure}

We have
$$\Jac=-a_2b_1xy-a_2b_2xy^3+2a_1b_2x^2y^2-2a_2b_3x^2y^3+2(a_1b_3-a_3b_1)\,x^3y^2+2a_3b_2x^3y^4,$$
where we can ignore the third and the fourth terms (see \rr{jac}).
Now the map \re{resmap} written in the above bases has the
following $15\times 17$ matrix, which we denote by $\mathbf A$.

\begin{equation}
 \left[
{\begin{array}{rrrrccccccccccccc} 0 & 0 & 0 & 0 & {a_{2}} & 0 & 0
& 0 & 0 & 0 & {b_{1}} & 0 & 0 & 0 & 0 & 0 &  - {a_{2}}{b_{1}} \\
0 & 0 & 0 & 0 & 0 & {a_{2}} & 0 & 0 & 0 & 0 & 0 & 0 & 0 & 0 & 0 & 0 & 0 \\
0 & 0 & 0 & 0 & 0 & 0 & {a_{2}} & 0 & 0 & 0 & {b_{2}} & 0 & 0 & 0 & 0 & 0 &  - {a_{2}}{b_{2}} \\
0 & 0 & 0 & 0 & {a_{1}} & 0 & 0 & 0 & 0 & 0 & 0 & {b_{1}} & 0 & 0 & 0 & 0 & 0 \\
1 & 0 & 0 & 0 & 0 & {a_{1}} & 0 & 0 & 0 & 0 & 0 & 0 & {b_{1}} & 0 & 0 & 0 & 0 \\
0 & 1 & 0 & 0 & 0 & 0 & {a_{1}} & {a_{2}} & 0 & 0 & 0 & {b_{2}} & 0 & 0 & 0 & 0 & 0 \\
0 & 0 & 1 & 0 & 0 & 0 & 0 & 0 & {a_{2}} & 0 & {b_{3}} & 0 & {b_{2}} & 0 & 0 & 0 &  0 \\
0 & 0 & 0 & 0 & 0 & 0 & 0 & {a_{1}} & 0 & 0 & 0 & 0 & 0 & {b_{1}} & 0 & 0 & 0 \\
0 & 0 & 0 & 0 & {a_{3}} & 0 & 0 & 0 & {a_{1}} & 0 & 0 & {b_{3}} & 0 & 0 & {b_{1}} & 0 &  2({a_{1}}{b_{3}}- {a_{3}}{b_{1}})  \\
0 & 0 & 0 & 1 & 0 & {a_{3}} & 0 & 0 & 0 & {a_{2}} & 0 & 0 & {b_{3}} & {b_{2}} & 0 & 0 & 0 \\
0 & 0 & 0 & 0 & 0 & 0 & {a_{3}} & 0 & 0 & 0 & 0 & 0 & 0 & 0 & {b_{2}} & 0 & 2{a_{3}}{b_{2}} \\
0 & 0 & 0 & 0 & 0 & 0 & 0 & 0 & 0 & {a_{1}} & 0 & 0 & 0 & 0 & 0 & {b_{1}} & 0 \\
0 & 0 & 0 & 0 & 0 & 0 & 0 & {a_{3}} & 0 & 0 & 0 & 0 & 0 & {b_{3}} & 0 & 0 & 0 \\
0 & 0 & 0 & 0 & 0 & 0 & 0 & 0 & {a_{3}} & 0 & 0 & 0 & 0 & 0 & {b_{3}} & {b_{2}} & 0 \\
0 & 0 & 0 & 0 & 0 & 0 & 0 & 0 & 0 & {a_{3}} & 0 & 0 & 0 & 0 & 0 &
 {b_{3}} & 0
\end{array}}
 \right]\nonumber
\end{equation}

It remains to solve the system {${\mathbf A}{\mathbf c}={\mathbf
b}$}, where $\mathbf c$ is the vector of unknowns (the
coefficients of the $h_i$ and $c$) and $\mathbf b$ is the monomial
$g=x^5y^4$ written in the basis for $S_{\tilde\D^\circ}$, i.e.
$\mathbf b=(0,\dots,0,1)^t$. With the help of {\tt Maple} we
obtain
$$c=\frac{1}{4}\,\frac{{a_{1}}^{2}\,{b_{2}}}{{a_{3}}\,( {a_{1}}\,{b_{3}} - {a_{3}}\,{b_{1}})
^{2}}.$$ Since the mixed volume (area) of $\D_1$ and $\D_2$ equals
4, we conclude that
$$\gres(x^5y^4)=\frac {{a_{1}}^{2}\,{b_{2}}}{{a_{3}}\,( {a_{1}}\,{b_{3}} - {a_{3}}\,{b_{1}})
^{2}}.$$
\end{Ex}

\section{Some geometry}\label{S:geometry}

The first step of our algorithm is constructing a lattice
$n$-dimensional polytope $\D_0$ such that the Minkowski sum
$\tilde\D=\D_0+\D$ contains both $\D$ and $\D(g)$ in its interior.
There are many ways of doing that. For example, one can take
$\D_0$ to be a sufficiently large dilate of $\D$ (translated so it
contains the origin in the interior). In general, this could
result in an unnecessarily large dimension of $S_{\tilde\D^\circ}$,
which determines the size of the linear system \re{linsys}.
Therefore, to minimize the size of the linear system one would
want to solve the following problem.

\begin{Prob} Given two lattice polytopes $\D$ and $\D'$ in $\R^n$, $\dim \D=n$,
find an $n$-dimensional lattice polytope $\D_0$ such that $\D+\D_0$ contains
both $\D$ and $\D'$ in its interior and has the smallest possible number of interior
lattice points.
\end{Prob}

This appears to be a hard optimization problem. Instead we will
consider a less challenging one. First, since the global residue
is linear we can assume that $g$ is a monomial, i.e. $\D(g)$ is a
point.

\begin{Prob}\label{Pr:segment} Given a convex polytope $\D$ and a point $u$ in $\R^n$
find a segment $I$ starting at the origin such that $u$ is
contained in the Minkowski sum $\D+I$ and the volume of $\D+I$ is
minimal.
\end{Prob}

If $I=[\hspace{1pt}0,m]$, for $m\in\Z^n$, is such a segment then we can take
$\D_0$ to be a ``narrow'' polytope with $0$ in the interior and
$m$ one of the vertices, as in \rf{ex2}. Then the volume (and
presumably the number of interior lattice points) of
$\tilde\D=\D+\D_0$ will be relatively small.

We will now show how \rpr{segment} can be reduced to a linear programming problem.
Let $\D\subset\R^n$ be an $n$-dimensional polytope, $I=[\hspace{1pt}0,v]$ a segment, $v\in\R^n$.
First, notice that the volume of $\D+I$ equals
$$\Vol_n(\D+I)=\Vol_n(\D)+|I|\cdot \Vol_{n-1}(\text{pr}_I\D),$$
where $\text{pr}_I\D$ is the projection of $\D$ onto the
hyperplane orthogonal to $I$, $\Vol_k$ the $k$-dimensional volume,
and $|I|$ the length of $I$. For each facet $\G\subset\D$ let
$n_\G$ denote the outer normal vector whose length equals the
($n-1$)-dimensional volume of $\G$. Then we can write
$$|I|\cdot\Vol_{n-1}(\text{pr}_I(\D))=\frac{1}{2}\sum_{\G\subset\D}|\langle n_\G,v\rangle|.$$
But the latter is the support function $h_Z$ of a convex polytope (zonotope) $Z$,
which is the Minkowski sum of segments:
$$h_Z(v)=\sum_{\G\subset\D}|\langle n_\G,v\rangle|,\quad Z=\sum_{\G\subset\D}[-n_\G,n_\G].$$
Indeed, $h_Z$ is the sum of the support functions of the segments. Also it is clear that
$$h_{[-n_\G,n_\G]}(v)=\max_{-1\leq t\leq 1}\langle\hspace{1pt} tn_\G,v\rangle=|\langle n_\G,v\rangle|.$$
The following figure shows the polytopes $\D$ and $Z$, and the normal fan $\Sig_Z$ of $Z$.

\begin{figure}[h]
\centerline{
\scalebox{0.55}{
\input{fig1.pstex_t}}}
\caption{}
\label{F:fig1}
\end{figure}


Now we get back to \rpr{segment}. After translating everything by $-u$ we may assume that
$u$ is at the origin. Then \rpr{segment} is equivalent to finding $x\in\D$ such that
the volume of $\D+[\hspace{1pt}0,-x]$ is minimal, which by the previous
discussion means minimizing the support function $h_Z(-x)=h_Z(x)$ on $\D$.

We can interpret this geometrically. The function $h_Z$ is a
nonnegative continuous function, linear on every cone of the normal
fan $\Sig_Z$. Its graph above the polytope $\D$ is a
``convex down'' polyhedral set in $\R^{n+1}$ (see \rf{fig2}). The set of points
with the smallest last coordinate is a face of this polyhedral
set. The projection of this face to $\D$ gives the solution to our
minimization problem.

\begin{figure}[h]
\centerline{
\scalebox{0.80}{
\input{fig2.pstex_t}}}
\caption{}
\label{F:fig2}
\end{figure}

Finally, note that the normal fan $\Sig_Z$ has a simple description.
It is obtained by translating all the facet hyperplanes $H_\G$,
for $\G\subset\D$, to the origin.


\end{document}